\documentclass[3p,12pt]{elsarticle}

\usepackage{fleqn}

\usepackage{geometry}
\usepackage{natbib}
\usepackage{pifont}
\usepackage{pifont}
\usepackage{natbib}
\usepackage{amssymb}

\usepackage[english]{babel}
\usepackage{pgf,pgfarrows,pgfnodes,pgfautomata,pgfheaps}
\usepackage{amsmath,amssymb}
\usepackage[latin1]{inputenc}
\usepackage{graphicx,amsmath,amsfonts,amssymb,amsthm,mathrsfs}
\input xy

\usepackage{color}
\usepackage[colorlinks]{hyperref}
\usepackage{amsmath}
\newtheorem{theorem}{Theorem}[section]

\newtheorem{definition}[theorem]{Definition}

\newtheorem{proposition}[theorem]{Proposition}

\newtheorem{remark}[theorem]{Remark}
\newtheorem{example}[theorem]{Example}

\numberwithin{equation}{section}

\journal{.  .  .  }
\begin{document}

\begin{frontmatter}

\title{ Controlled weaving frames in Hilbert spaces
}


\author{R. Rezapour \footnote{ Department of Mathematics, Shabestar Branch, Islamic Azad University, Shaberstar, Iran.,  \ Email: {\tt  rrezapour@marandiau.ac.ir}} ,
A. Rahimi \footnote{Corresponding outhor. Department of Mathematics, University of Maragheh, Maragheh, Iran.,  \ Email: {\tt rahimi@maragheh.ac.ir}},
E. Osgooei\footnote{Department of sciences, Urmia University of Technology, Urmia, Iran \ Email: {\tt e.osgooei@uut.ac.ir }}
 ,
H. Dehghan \footnote{Department of Mathematics, Institute for Advanced Studies in Basic Sciences (IASBS), Gava Zang, Zanjan 45137-66731, Iran, \ Email: {\tt hossein.dehgan@gmail.com }}}

%
%
%
%
%
%



\begin{abstract}

In this paper, we first introduce the notion of controlled weaving K-g-frames in Hilbert spaces. Then, we present sufficient conditions for controlled weaving K-g-frames in separable Hilbert spaces. Also, a characterization of controlled weaving K-g-frames is given in terms of an operator.  Finally,  we show that if bounds of frames associated with atomic spaces are positively confined, then controlled K-g-woven frames gives ordinary weaving K-frames and vice-versa.
\end{abstract}


\end{frontmatter}



\section{Introduction and Preliminaries}

The Hilbert space is the natural framework for the mathematical description of many
areas of physics: certainly for quantum mechanics (and quantum field theory), and
also for classical electrodynamics (in the energy norm), classical scattering theory \cite{Lax},
signal and image analysis (finite energy signals) \cite{Lynn}, etc. In each cases, there arises the
problem of representing an arbitrary vector in terms of simpler ones, i.e., in terms of
the elements of some basis $\{f_j\}_{j\in \mathbb{N}}$. The most economical solution, and the one
advocated by mathematicians, is of course to use an orthonormal basis, which gives in
addition the uniqueness of the decomposition,
\begin{eqnarray*}
f= \sum_{j\in \mathbb{N}} \langle f_j, f \rangle f_j
\end{eqnarray*}
for any $f$ in the underlying Hilbert space.
Unfortunately, orthonormal bases are often difficult to find and sometimes hard to work
with. One way to give up orthogonality of the basis vectors and uniqueness
of the decomposition- while maintaining its other useful properties: fast convergence,
numerical stability of the reconstruction $ \langle f_j, f \rangle \to f$ and etc- is using the notion of frames. The notion of frames in Hilbert spaces was introduced  by Duffin and Schaeffer during their study of nonharmonic Fourier series in 1952 \cite{Duffin}.
Let $ \mathcal{H}$ be a separable Hilbert space with the inner product $\langle . , .\rangle$. Recall that a countable family of elements $\{f_j\}_{j\in J}$ in $\mathcal{H}$  is a frame for $\mathcal{H}$ if there exist constants $A, B > 0$ such that
\begin{eqnarray}\label{111}
A \Vert f\Vert^{2}\leq\sum_{j\in J} |\langle f_j, f\rangle|^2 \leq B \Vert f\Vert^{2},\ \  \forall f\in \mathcal{H}.
\end{eqnarray}
The constants $A$ and $B$ are called lower and upper frame bounds, respectively. In case $A=B$, it called a tight frame and if $A=B=1$ it is known Parseval frame. If the second inequality in (\ref{111}) holds, it called a Bessel sequence. For a frame $\{f_j\}_{j\in J}$ in $\mathcal{H}$, the operator $Sf = \sum_{j\in J} \langle f_j, f\rangle f_j, f\in\mathcal{H} $ called the frame operator. This operator is a bounded, self-adjoint, invertible and positive operator and any $f\in\mathcal{H}$ has an expansion
\begin{eqnarray}\label{222}
f = \sum_{j\in J} \langle S^{-1} f_j, f\rangle f_j= \sum_{j\in J} \langle  f_j, f\rangle S^{-1}f_j.
\end{eqnarray} The family $\{S^{-1}f_j\}_{j\in J}$ is also a frame with bounds $B^{-1}, A^{-1}$, this frame is called the canonical dual or reciprocal
frame of $\{f_j\}_{j\in J}$.
\par
After the fundamental paper by Daubechies, Grossmann and Meyer (2017 Abel Prize winer)  in 1986 \cite{Daubechies}, the theory of frames found more attentions and a lot of papers and books published in this area. They took the key step of connecting frames with wavelets and Gabor systems in that paper.
\par
In practice, the expansion (\ref{222}) converges quite fast and can be truncated after a few
terms, provided $|B/A - 1| \ll 1$. This fact is the key to all the reconstruction formulae, used in signal or image processing in terms of wavelet or Gabor analysis. For more details we refer to \cite{Daubechies} and also to the review papers \cite{Heil,Heil2}. For a more complete treatment of frame theory we recommend the excellent book of
Christensen \cite{Christensen}, and the tutorials of Casazza \cite{Ca1,Ca2} and the Memoir of Han
and Larson \cite{han1}.

\par

We denote by $\mathcal{L(H)}$  the set of all bounded linear operators. A bounded operator $T\in \mathcal{L(H)}$ is called positive (respectively, non- negative), if $\langle Tf, f\rangle>0$ for all $f\neq 0$ (respectively, $\langle Tf, f\rangle\geq 0$ for all $f$ ). Every non-negative operator is clearly self-adjoint.   If $T\in \mathcal{L(H)}$ is non-negative then there exists a unique non-negative operator $S$ such that $S^{2}=T$. This will be denoted by $S=T^{\frac{1}{2}}$. Moreover, if an operator $D$ commutes with $T$ then $D$ commutes with every operator in the $C^*$-algebra generated by $T$ and $I$, specially $D$ commutes with  $T^{\frac{1}{2}}$. Let $\mathcal{L^{+}(H)}$ be the set of positive operators on $H$. For self-adjoint operators $T_{1}$ and $T_{2}$, the notation $T_{1}\leq T_{2}$ or $T_{2}-T_{1}\geq 0$ means
$$
\langle T_{1}f, f\rangle\leq\langle T_{2}f, f\rangle,  \quad f\in \mathcal{H}.
$$
We denote by $\mathcal{GL(H)}$  the set of all bounded linear operators which have
bounded inverse.
It is easy to see that if $S, T\in \mathcal{GL(H)}$ then $T^{*}, T^{-1}$ and $ST$ are also in $\mathcal{GL(H)}$. Let $\mathcal{GL^{+}(H)}$ be the set of all positive operators in $\mathcal{GL(H)}$. For $u\in \mathcal{L(H)}$, $u\in \mathcal{GL^{+}(H)}$ if and only if there exists positive constants $0<m \leq M<\infty$ such that
$$mI\leq u \leq MI.$$ For $u^{-1}$, $$M^{-1}I \leq u^{-1} \leq m^{-1}I.$$

Throughout this paper,  $\mathcal{K}_{1}$, $\mathcal{K}_{2}$ and $\mathcal{H}$ are complex separable Hilbert spaces, $K\in \mathcal{L(H)}$,  $C, C'\in \mathcal{GL^{+}(H)}$ and  $\{\mathcal{H}_{j}\}_{j=1}^{\infty}\subset \mathcal{K}_{1}$ and $\{\mathcal{W}_{k}\}_{k=1}^{\infty}\subset \mathcal{K}_{2}$ are sequences of closed subspaces.


The following theorem can be found in \cite{Murphy}.
\begin{theorem}\label{ineqa1}
 Let $T_1, T_2, T_3 \in \mathcal{L(H)}$ and $T_1 \leq T_2$. Suppose $T_3 \geq 0$ commutes with $T_1$ and $T_2$ then
$$T_1T_3\leq T_2T_3.$$
\end{theorem}

\par

In the last decade, motivated by new applications of frame theory, many generalizations of
frames introduced: fusion frames, continuous frames, Banach frames, g-frames, K-frames, operator valued frames, p-frames, pg-frames, frames for Hilbert $C^*$-modules and etc.
\subsection{K-frames}
  Atomic systems for subspaces
were first introduced by Feichtinger and Werther in \cite{Feichtinger} based on examples
arising in sampling theory.
  In 2011, Gavruta \cite{Gavruta} introduced K-frames in
Hilbert spaces to study atomic decomposition systems, and discussed some
properties of them. Let $K$ be a linear and bounded operator on $ \mathcal{H}$.
A family of elements $\{f_j\}_{j\in J}$ in $\mathcal{H}$  is a $K$-frame for $\mathcal{H}$ if there exist constants $A, B > 0$ such that
\begin{eqnarray}\label{Gframe}
A \Vert K^* f\Vert^{2}\leq\sum_{j\in J} |\langle f_j, f\rangle|^2 \leq B \Vert f\Vert^{2},\ \  \forall f\in \mathcal{H}.
\end{eqnarray}
  K-frames are limited to the range of a bounded linear
operator in Hilbert spaces, that is, they replace the lower bound condition
$A\Vert f\Vert^2$ of classical frames (\ref{Gframe}) by new lower condition $A\Vert K^* f\Vert^2$.   In recent years,
K-frames have been widely studied in \cite{Asgari,Gavruta2,Xiao} with a paramount field
in frame theory.
\subsection{Controlled frames}
Controlled frames, as one of the newest generalizations of
frames, have been introduced to improve the numerical efficiency of iterative
algorithms for inverting the frame operator on abstract Hilbert spaces \cite{Balazs},
however, they are used earlier just as a tool for spherical wavelets \cite{Bogdanova}. Since then, controlled frames have been widely studied.
In 2016, Hua and Huang \cite{Hua} introduced $(C, C')$-controlled K-g-frame as follows: Assume that $K, C, C'$ be  linear and bounded operators on $ \mathcal{H}$ such that $C$ and $C'$ are positive and have bounded inverse.
 The family $\{\Lambda_{j} : \mathcal{H} \to \mathcal{K}_{j} \}_{j\in J}$ is called $(C, C')$-controlled K-g-frame for $\mathcal{H}$ with respect to $\{\mathcal{K}_{j}\}_{j\in J}$, if there exist constants $ 0<A\leq B<\infty$ such that
\begin{eqnarray*}
A  \Vert K^{*}f\Vert^{2}\leq\sum_{j\in J}\langle \Lambda_{j}Cf, \Lambda_{j}C'f\rangle\leq B \Vert f\Vert^{2},\ \  \forall f\in \mathcal{H}.
\end{eqnarray*}
\subsection{Weaving frames}
Very recently, Bemrose et al. \cite{Bemrose} introduced a new concept of "weaving frames" in
separable Hilbert spaces which is motivated by the following question in distributed signal processing: given are two sets $\{f_j\}_{j\in J}$ and $\{g_j\}_{j\in J}$ of linear measurements with stable recovery guarantees,
in mathematical terminology each set is a frame labeled by a node or sensor $j\in J$.
At each sensor a signal $f$ was measured either with $f_j$ or with $g_j$, so that the collected
information is the set of numbers $\{ \langle f, f_j \rangle\}_{j\in \sigma}\cup \{ \langle f, g_j \rangle\}_{j\in \sigma^c}$ for some subset $\sigma \subseteq J$. Can $f$ still be recovered robustly from these measurements, no matter which kind of measurement
has been made at each node? In other words, is the set $\{  f_j \}_{j\in \sigma}\cup \{ g_j \}_{j\in \sigma^c}$ a
frame for all subsets $\sigma \subseteq J$? This question led them to define woven frame as:
Two frames $\{f_j\}_{j\in J}$ and $\{g_j\}_{j\in J}$ for a Hilbert space $\mathcal{H}$ are (weakly) woven if for every subset
$\sigma \subseteq J$, the family $\{  f_j \}_{j\in \sigma}\cup \{ g_j \}_{j\in \sigma^c}$ is a frame for $\mathcal{H}$.
Notions related to weaving frames, such as weaving Riesz bases, weaving frames in Banach space setting, weaving vector-valued frames and continuous weaving frames  were considered by many researchers. For more details we refer the reader to \cite{ Casazza, Casazza2, Daubechies, Vashisht} and references there in.
\par
In this paper, motivated and inspired by the above mentioned works we introduce the concept of controlled weaving K-g-frame. This frame includes ordinary frame, $K$-frame, $g$-frame, controlled frame and weaving frame. We present characterization theorems of controlled weaving K-g-frames and construt an example to illustrate our results. We provide a necessary and sufficient condition for $(C, C')$-controlled $K$-g-woven frames which connects to ordinary weaving $K$-frames.  Through investigating the characterizations
of $(C,C')$-controlled weaving K-g-frames, we obtain some equivalent conditions
of $(C,C')$-controlled weaving K-g-frames.
\par


\begin{definition}\label{D-Hua}
(\cite{Hua}) Assume that $K\in \mathcal{L(H)}$ and $C, C'\in \mathcal{GL^{+}(H)}$.   The family $\Lambda=\{\Lambda_{j}\}_{j\in J}$ is called $(C, C')$-controlled K-g-frame for $\mathcal{H}$ with respect to $\{\mathcal{H}_{j}\}_{j\in J}$, if there exist constants $ 0<A\leq B<\infty$ such that
\begin{eqnarray}\label{CCKgf}
A  \Vert K^{*}f\Vert^{2}\leq\sum_{j\in J}\langle \Lambda_{j}Cf, \Lambda_{j}C'f\rangle\leq B \Vert f\Vert^{2}, \forall f\in \mathcal{H}.
\end{eqnarray}
The constants $A$ and $B$ are called lower and upper frame bounds for $(C, C')$-controlled K-g-frame, respectively.
\end{definition}
If the right-hand side of (\ref{CCKgf}) holds then we call $\Lambda=\{\Lambda_{j}\}_{j\in J}$ a  $(C, C')$-controlled $K$-g-Bessel sequence for $\mathcal{H}$ with respect to$\{\mathcal{H}_{j}\}_{j\in J}$.  \\
If $C'=I$ then we call $\Lambda=\{\Lambda_{j}\}_{j\in J}$ a   $C$-controlled K-g-frame for $\mathcal{H}$ with respect to
$\{\mathcal{H}_{j}\}_{j\in J}$.
If $K=C=C'=I$ in (\ref{CCKgf}) then $(C, C')$-controlled K-g-frame will be the g-frame.  \\

Associated with a $(C, C')$-controlled $g$-Bessel sequence $\Lambda$, we shall denote the representation space as follows:
$$
\left(\sum_{n\in \mathbb{N}}\oplus \mathcal{H}_{n}\right)_{\ell^{2}}=\left\{\{g_{n}\}_{n=1}^{\infty}:g_{n}\in \mathcal{H}_{n}(n\in \mathbb{N}),\sum_{n=1}^{\infty}\Vert g_{n}\Vert^{2}<+\infty\right\}.
$$
The operator $T_{(C,C')}$ : $\left( \sum_{n\in \mathbb{N}}\oplus \mathcal{H}_{n}\right)_{\ell^{2}}\rightarrow \mathcal{H}$ defined by
$$
T_{(C,C')}\left(\{g_{n}\}_{n=1}^{\infty}\right)=\sum_{n=1}^{\infty}(CC')^{\frac{1}{2}} \Lambda_{n}^{*}g_{n}
$$
is called the  pre-frame operator (or  synthesis operator) and the  analysis operator
$$
T^{*}_{(C,C')}:\mathcal{H}\rightarrow \left(\sum_{n\in \mathbb{N}}\oplus \mathcal{H}_{n}\right)_{\ell^{2}},
$$
which is the adjoint of $T_{\Lambda}$, is defined by
$$
T^{*}_{(C,C')} : f\rightarrow\{\Lambda_{n} (CC')^{\frac{1}{2}} f\}_{n=1}^{\infty},\ f\in \mathcal{H}.
$$
The {\it frame operator} $S_{\Lambda}$ associated with $\Lambda$ is defined as
$$
S_{(C,C')}=T_{(C,C')}T^{*}_{(C,C')}:\mathcal{H}\rightarrow \mathcal{H}
$$
$$
S_{(C,C')}:f\rightarrow\sum_{n=1}^{\infty}C'\Lambda_{n}^{*}\Lambda_{n}Cf,\ f\in \mathcal{H}.
$$
If $\Lambda$ is a frame for $\mathcal{H}$, then $S_{(C,C')}$ is a linear, bounded and positive operator.

Throughout this paper, the mapping $S_\Lambda : \mathcal{H} \to \mathcal{H}$ defined by
$S_\Lambda f  :=\sum_{j=1}^\infty \Lambda^*_{j} \Lambda_j f $ is well-defined and
  commutes with $C$  and $C^\prime$.

\begin{definition}
Let $\{\mathcal{H}_{n}\}_{n\in \mathbb{N}}$ be a sequence of closed subspaces of $\mathcal{H}$. The family of bounded linear operators $\{\Xi_{n}\}_{n=1}^{\infty}$ defined by
$$
\Xi_{n}:\left(\sum_{k=1}^{\infty}\oplus \mathcal{H}_{k}\right)_{\ell^{2}}\rightarrow \mathcal{H},\ \ \ \ \ \ \Xi_{n}\left(\{g_{k}\}_{k=1}^{\infty}\right)=g_{n}
$$
is called the standard $g$-orthonormal basis for $ \left(\sum_{n=1}^{\infty}\oplus \mathcal{H}_{n}\right)_{\ell^{2}}$ with respect to $\{\mathcal{H}_{n}\}_{n\in \mathbb{N}}.$

\end{definition}

\begin{remark}
If $\{\Xi_{n}\}_{n=1}^{\infty}$ is the standard $g$-orthonormal basis for $\left( \sum_{n=1}^{\infty}\oplus \mathcal{H}_{n}\right)_{\ell^{2}}$, then
$$\Xi_{j}\Xi_{i}^{*}=\left\{\begin{array}{ll}
I_{o},\ \ \  i=j, & \\
  &  \\
0,\ \ \   i\neq j & \end{array}\right.\ \  \  \mbox{and} \ \ \ \ \ \ \sum_{j\in \mathbb{N}}\Xi_{j}^{*}\Xi_{j}=I_{\ell^{2}},
$$
where $I_{\ell^{2}}$ denotes the identity operator on $\left( \sum_{n=1}^{\infty}\oplus \mathcal{H}_{n}\right)_{\ell^{2}}$ and $I_{o}$ denotes the identity operator on $\mathcal{H}_{j}.$
\end{remark}


\section{Main Results}

We firs define the weaving $(C, C')$-controlled  $K$-$g$-frames in separable Hilbert spaces.

\begin{definition}\label{D-WCKgf}
 Two $(C, C')$-controlled $K$-$g$-frames $\Lambda=\{\Lambda_{j}\}_{j=1}^{\infty}$ and $\Omega=\{\Omega_{k}\}_{k=1}^{\infty}$ for $\mathcal{H}$ with respect to $\{\mathcal{H}_{j}\}_{j=1}^{\infty}$ and $\{\mathcal{W}_{k}\}_{k=1}^{\infty}$ (respectively) are $(C, C')$-controlled $K$-$g$-woven if there are universal constants $ 0<A\leq B<\infty$ so that for every subset $\sigma$ of $\mathbb{N}$, the family $\{\Lambda_{j}\}_{j\in\sigma}\cup\{\Omega_{k}\}_{k\in\sigma^{c}}$ is a $(C, C')$-controlled $K$-$g$-frame for $\mathcal{H}$ (with respect to $\{\mathcal{H}_{j}\}_{j\in\sigma}\cup\{\mathcal{W}_{k}\}_{k\in\sigma^{c}}$) with lower and upper $K$-$g$-frame bounds $A$ and $B$, respectively.
\end{definition}
The following  proposition show that every family of $(C, C')$-controlled $K$-$g$-frames has a universal $(C, C')$-controlled $K$-$g$-Bessel bound.

\begin{proposition}\label{3.2}
Let $\Lambda=\{\Lambda_{j}\}_{j=1}^{\infty}$   and $\Omega=\{\Omega_{k}\}_{k=1}^{\infty}$   be $(C, C')$-controlled $K$-$g$-Bessel sequences in $\mathcal{H}$ (with  respect to $\{\mathcal{H}_{j}\}_{j=1}^{\infty}$ and $\{\mathcal{W}_{k}\}_{k=1}^{\infty},$   respectively)   with bound $B_{1}$   and $B_{2}$,   respectively. Then, for any subset $\sigma$   of $\mathbb{N}$,   the family $\{\Lambda_{j}\}_{j\in\sigma}\cup\{\Omega_{k}\}_{k\in\sigma^{c}}$   is a $(C, C')$-controlled $K$-$g$-Bessel sequence with bound $B_{1}+B_{2}.$
\end{proposition}

\noindent\textbf{Proof.}
Let $\sigma$ be an arbitrary subset of $\mathbb{N}$. For each $f\in \mathcal{H}$, we have
\begin{align*}
\sum_{j\in\sigma}\langle \Lambda_{j}Cf, \Lambda_{j}C'f\rangle +\sum_{k\in\sigma^{c}}\langle \Omega_{k}Cf, \Omega_{k}C'f\rangle &\leq\sum_{j\in \mathbb{N}} \langle \Lambda_{j}Cf, \Lambda_{j}C'f\rangle+\sum_{k\in \mathbb{N}}\langle \Omega_{k}Cf, \Omega_{k}C'f\rangle\\
& \leq(B_{1}+B_{2})\Vert f\Vert^{2}.
\end{align*}
This is the required universal upper $(C, C')$-controlled $K$-$g$-frame bound for $\Lambda$ and $\Omega$.
$\Box$\\\\
The following theorem is a characterization of $(C, C')$-controlled $K$-$g$-woven frames in terms of an operator.

\begin{theorem}
Let $\Lambda=\{\Lambda_{j}\}_{j=1}^{\infty}$  and $\Omega=\{\Omega_{j}\}_{j=1}^{\infty}$  be sequences of operators such that $\Lambda_{j}\in B(\mathcal{H},\ \mathcal{H}_{j})$  and $\Omega_{j}\in B(\mathcal{H},\ \mathcal{W}_{j})$  for all $j\in \mathbb{N}$.  Then the following conditions are equivalent.
\begin{itemize}
\item[(i)]
$\Lambda$  and $\Omega$  are $(C, C')$-controlled $K$-$g$-woven frames for $\mathcal{H}.$
\item[(ii)]
\begin{itemize}
\item[(a)]
There exists $A>0$  such that for any subset $\sigma$  of $\mathbb{N}$  there exists a bounded linear operator $U_{\sigma}$ : $\left(\displaystyle \sum_{n=1}^{\infty}\oplus \mathcal{Z}_{n}^{\sigma}\right)_{\ell^{2}}\rightarrow \mathcal{H}$ ( here, $\mathcal{Z}_{n}^{\sigma}=\mathcal{H}_{n}$  for $ n\in\sigma$  and $\mathcal{Z}_{n}^{\sigma}=\mathcal{W}_{n}$  for $n\in\sigma^{c}$)  such that
\begin{align*}
U_{\sigma}\left(\Xi_{j}^{*}\Xi_{j}\{g_{n}\}_{n=1}^{\infty}\right)=\left\{\begin{array}{l}
(CC')^{\frac{1}{2}}\Lambda_{j}^{*}(g_{j})\ ,\ j\in\sigma\\
   \\
(CC')^{\frac{1}{2}}\Omega_{j}^{*}(g_{j})\ ,\ j\in\sigma^{c}
\end{array}\right.
\end{align*}
 for all
\begin{align*}
\{g_{n}\}_{n=1}^{\infty} \in\left(\sum_{n=1}^{\infty}\oplus \mathcal{Z}_{n}^{\sigma}\right)_{\ell^{2}},
\end{align*}
 where $\{\Xi_{n}\}_{n=1}^{\infty}$  is the  standard $g$-orthonormal basis for
$$
\left(\sum_{n=1}^{\infty}\oplus \mathcal{Z}_{n}^{\sigma}\right)_{\ell^{2}}.
$$
\item[(b)]
$A K K^{*}\leq U_{\sigma}U_{\sigma}^{*}.$
\end{itemize}
\end{itemize}
\end{theorem}
\noindent\textbf{Proof.}
 $(i)\Rightarrow(ii)$ : Let $A$ be the universal lower $(C, C')$-controlled $K$-$g$-frame bound for $\Lambda$ and $\Omega$. For each subset $\sigma$ of $\mathbb{N}$, let $T_{\sigma}$ be the pre-frame operator associated with $\{\Lambda_{j}\}_{j\in\sigma}\cup\{\Omega_{j}\}_{j\in\sigma^{c}}.$
Put $U_{\sigma}=T_{\sigma}$. Then,
\begin{align*}
U_{\sigma}(\Xi_{j}^{*}\Xi_{j}\{g_{n}\}_{n=1}^{\infty}) & =T_{\sigma}(\Xi_{j}^{*}\Xi_{j}\{g_{n}\}_{n=1}^{\infty}) \\
&=\left\{\begin{array}{l}
(CC')^{\frac{1}{2}}\Lambda_{j}^{*}(g_{j})\  ,\ j\in\sigma   \\
(CC')^{\frac{1}{2}}\Omega_{j}^{*}(g_{j})\  ,\ j \in \sigma^{c}
\end{array}\right. , \ \ \ \ \ \ \{g_{n}\}_{n=1}^{\infty}\in\left(\sum_{n=1}^{\infty}\oplus \mathcal{Z}_{n}^{\sigma}\right)_{\ell^{2}}.
\end{align*}
Moreover, for all $f\in \mathcal{H}$ we have
\begin{align*}
A \langle K K^{*}f,\ f\rangle  & =  A\Vert K^{*}f\Vert^{2}\leq
\sum_{j\in\sigma}\langle \Lambda_{j}Cf, \Lambda_{j}C'f\rangle +\sum_{k\in\sigma^{c}}\langle \Omega_{k}Cf, \Omega_{k}C'f\rangle \\
& = \sum_{j\in\sigma}   \left\langle  f, C\Lambda_{j}^*\Lambda_{j}C' f\right\rangle  +\sum_{k\in\sigma^{c}}\langle f, C \Omega_{k}^*\Omega_{k}C'f\rangle \\
& = \sum_{j\in\sigma}   \left\langle \Lambda_{j}(CC')^{\frac{1}{2}} f, \Lambda_{j}(CC')^{\frac{1}{2}} f\right\rangle +\sum_{k\in\sigma^{c}}  \left\langle \Omega_{j}(CC')^{\frac{1}{2}} f, \Omega_{j}(CC')^{\frac{1}{2}} f\right\rangle\\
& = \sum_{j\in\sigma}   \Vert\Lambda_{j} (CC')^{\frac{1}{2}} f \Vert^{2} +\sum_{k\in\sigma^{c}}  \Vert\Omega_{k} (CC')^{\frac{1}{2}} f \Vert^{2} \\
&=\Vert T_{\sigma}^{*}(f)\Vert^{2}=\Vert U_{\sigma}^{*}f\Vert^{2}.
\end{align*}
This gives $A K K^{*}\leq U_{\sigma}U_{\sigma}^{*}.$

\noindent $(ii)\Rightarrow(i)$ : Let $\sigma$ be an arbitrary subset of $\mathbb{N}$. First we find $U_{\sigma}^{*}$ and for this (by using $(a)$ ) we have
\begin{align*}
\left\langle U_{\sigma}(\{g_{n}\}_{n=1}^{\infty}),\ f\right\rangle & =\left\langle U_{\sigma}\left(\sum_{j=1}^{\infty}\Xi_{j}^{*}\Xi_{j}(\{g_{n}\}_{n=1}^{\infty})\right),\ f\right\rangle\\
&=\left\langle\sum_{j=1}^{\infty}U_{\sigma}(\Xi_{j}^{*}\Xi_{j}(\{g_{n}\}_{n=1}^{\infty})),\ f\right\rangle\\
&=\left\langle\sum_{j\in\sigma}(CC')^{\frac{1}{2}}\Lambda_{j}^{*}(g_{j})+\sum_{j\in\sigma^{c}}(CC')^{\frac{1}{2}}\Omega_{j}^{*}(g_{j}),\ f\right\rangle\\
&=\sum_{j\in\sigma}\left\langle g_{j},\ \Lambda_{j}(CC')^{\frac{1}{2}}(f)\right\rangle+\sum_{j\in\sigma^{c}}\left\langle g_{j},\ \Omega_{j}(CC')^{\frac{1}{2}}(f)\right\rangle\\
&=\left\langle\{g_{n}\}_{n\in \mathbb{N}},\ \{\Lambda_{i}(CC')^{\frac{1}{2}}(f)\}_{j\in\sigma}\cup\{\Omega_{i}(CC')^{\frac{1}{2}}(f)\}_{j\in\sigma^{c}}\right\rangle
\end{align*}
for all $f\in \mathcal{H}, \displaystyle \{g_{n}\}_{n=1}^{\infty}\in\left(\sum_{n=1}^{\infty}\oplus \mathcal{Z}_{n}^{\sigma}\right)_{\ell^{2}}.$
This implies that
\begin{align}\label{Hd2}
U_{\sigma}^{*}(f)= & \left\{\Lambda_{j}(CC')^{\frac{1}{2}}f\right\}_{j\in\sigma}\cup\left\{\Omega_{j}(CC')^{\frac{1}{2}}f\right\}_{j\in\sigma^{c}}
\end{align}
 for all $f\in \mathcal{H}$.
By using ({\it b}) and (\ref{Hd2}), we have
\begin{align*}
A\Vert K^{*}f\Vert^{2} &= A\langle K K^{*}f,\ f\rangle\leq\langle U_{\sigma}U_{\sigma}^{*}f,\ f\rangle=\Vert U_{\sigma}^{*}f\Vert^{2}\\
& = \sum_{j\in\sigma}\Vert\Lambda_{j}(CC')^{\frac{1}{2}}f\Vert^{2}+\sum_{j\in\sigma^{c}}\Vert\Omega_{j}(CC')^{\frac{1}{2}}f\Vert^{2}
\end{align*}
 for all $f\in \mathcal{H}$.
This gives a universal lower $(C, C')$-controlled $K$-$g$-frame bound for the family $\Lambda$ and $\Omega.$

Next we prove that $\Lambda$ satisfies upper $(C, C')$-controlled $K$-$g$-frame inequality. For each $f\in \mathcal{H},$ we have
$$
\sum_{j\in \mathbb{N}}\Vert\Lambda_{j}(CC')^{\frac{1}{2}}f\Vert^{2}=\Vert U_{\mathbb{N}}^{*}f\Vert^{2}\leq\Vert U_{\mathbb{N}}^{*}\Vert^{2}\Vert f\Vert^{2}.
$$
Similarly, $\Omega$ satisfies upper $(C, C')$-controlled $K$-$g$-frame inequality. Thus, by applying Proposition \ref{3.2} we can obtain a universal upper $(C, C')$-controlled $K$-$g$-frame bound. Hence $\Lambda$ and $\Omega$ are $(C, C')$-controlled $K$-$g$-woven. $\square $

\begin{remark}
For each $\sigma\subseteq \mathbb{N}$, if $U_{\sigma}$ : $\displaystyle \left(\sum_{n=1}^{\infty}\oplus \mathcal{Z}_{n}^{\sigma}\right)_{\ell^{2}}\rightarrow \mathcal{H}$ (here, $\mathcal{Z}_{n}^{\sigma}=\mathcal{H}_{n}$ for $ n\in\sigma$ and $\mathcal{Z}_{n}^{\sigma}=\mathcal{W}_{n}$ for $n\in\sigma^{c}$) is a bounded linear operator such that

\begin{align*}
U_{\sigma}\left(\Xi_{j}^{*}\Xi_{j}\{g_{n}\}_{n=1}^{\infty}\right)=\left\{\begin{array}{l}
(CC')^{\frac{1}{2}}\Lambda_{j}^{*}(g_{j})\ ,\ j\in\sigma\\
   \\
(CC')^{\frac{1}{2}}\Omega_{j}^{*}(g_{j})\ ,\ j\in\sigma^{c}
\end{array}\right.
\end{align*}
 for all
\begin{align*}
\{g_{n}\}_{n=1}^{\infty} \in\left(\sum_{n=1}^{\infty}\oplus \mathcal{Z}_{n}^{\sigma}\right)_{\ell^{2}},
\end{align*}
Then, $U_{\sigma}^{*}(f)=\{\Lambda_{j}(CC')^{\frac{1}{2}}f\}_{j\in\sigma}\cup\{\Omega_{j}(CC')^{\frac{1}{2}}f\}_{j\in\sigma^{c}}$ for all $f\in \mathcal{H}.$

\end{remark}

\begin{theorem}\label{T3.6}
 Let  $\Lambda\equiv\{\Lambda_{j}\}_{j=1}^{\infty}$  and $\Omega\equiv\{\Omega_{j}\}_{j=1}^{\infty}$   be  $(C, C')$-controlled  $K$-$g$-frames for  $\mathcal{H}$   with respect to  $\{\mathcal{H}_{j}\}_{j=1}^{\infty}$   and with  $(C, C')$-controlled  $K$-$g$-frame bounds  $A, B$   and  $\alpha, \beta$,  respectively. Let  $\{e_{k}\}_{k=1}^{\infty}$   be an orthonormal basis for  $\mathcal{H}$   such that for each  $j\in \mathbb{N}, \{\Lambda_{j}(CC')^{\frac{1}{2}}(e_{k})\}_{k=1}^{\infty}$    and $\{\Omega_{j}(CC')^{\frac{1}{2}}(e_{k})\}_{k=1}^{\infty}$ are orthogonal sets. If for each  $k\in \mathbb{N}$,  there exists a sequence of real scalars  $\{\beta_{ij}^{k}\}_{i,j\in \mathbb{N}}$   such that
$$ \inf\{|\beta_{kj}^{k}|^{2} : j\in \mathbb{N}\}\geq M>0\ \ \ \  \mbox{and}\ \ \ \   \Omega_{j}(CC')^{\frac{1}{2}}(e_{k})=\sum_{i=1}^{\infty}\beta_{ij}^{k}\Lambda_{j}(CC')^{\frac{1}{2}}(e_{i}) (j\in \mathbb{N}),$$
then  $\Lambda$    and  $\Omega$  are  $(C, C')$-controlled $K$-$g$-woven with universal  $(C, C')$-controlled $K$-$g$-frame bounds  $(\min \{ 1,\ M\})$   and  $(B+\beta)$.
\end{theorem}
\noindent\textbf{Proof.}
 For each subset $\sigma\subset \mathbb{N}$ and for every $f=\displaystyle \sum_{k=1}^{\infty}\alpha_{k}e_{k}\in \mathcal{H}$, we have
\begin{align*}
(B+\beta)\Vert f\Vert^{2} &\geq\sum_{j\in\sigma}\langle \Lambda_{j}Cf, \Lambda_{j}C'f\rangle +
 \sum_{j\in\sigma^{c}}\langle \Omega_{j}Cf, \Omega_{j}C'f\rangle\\
& = \sum_{j\in\sigma}\langle \Lambda_{j}Cf, \Lambda_{j}C'f\rangle +
\sum_{j\in\sigma^{c}}\langle \Omega_{j}(CC')^{\frac{1}{2}}f, \Omega_{j}(CC')^{\frac{1}{2}}f\rangle\\
& =\sum_{j\in\sigma}\langle \Lambda_{j}Cf, \Lambda_{j}C'f\rangle
+
\sum_{j\in\sigma^{c}}\left\langle \Omega_{j}(CC')^{\frac{1}{2}}\left(\sum_{k=1}^{\infty}\alpha_{k}e_{k}\right), \Omega_{j}(CC')^{\frac{1}{2}}\left(\sum_{k=1}^{\infty}\alpha_{k}e_{k}\right)\right\rangle\\
& =\sum_{j\in\sigma}\langle \Lambda_{j}Cf, \Lambda_{j}C'f\rangle
+
\sum_{j\in\sigma^{c}}\left\langle \sum_{k=1}^{\infty}\alpha_{k} \Omega_{j}(CC')^{\frac{1}{2}}\left(e_{k}\right), \sum_{k=1}^{\infty}\alpha_{k} \Omega_{j}(CC')^{\frac{1}{2}}\left(e_{k}\right)\right\rangle\\
& =\sum_{j\in\sigma}\langle \Lambda_{j}Cf, \Lambda_{j}C'f\rangle
+
\sum_{j\in\sigma^{c}}\left\langle \sum_{k=1}^{\infty}\alpha_{k} \sum_{i=1}^{\infty}\beta_{ij}^{k}\Lambda_{j}(CC')^{\frac{1}{2}}(e_{i}), \sum_{k=1}^{\infty}\alpha_{k}\sum_{i=1}^{\infty}\beta_{ij}^{k}\Lambda_{j}(CC')^{\frac{1}{2}}(e_{i})\right\rangle\\
& =\sum_{j\in\sigma}\langle \Lambda_{j}Cf, \Lambda_{j}C'f\rangle
+
\sum_{j\in\sigma^{c}}  \sum_{k=1}^{\infty}|\alpha_{k}|^2 \sum_{i=1}^{\infty}|\beta_{ij}^{k}|^2 \left\langle  \Lambda_{j}(CC')^{\frac{1}{2}}(e_{i}), \Lambda_{j}(CC')^{\frac{1}{2}}(e_{i})\right\rangle\\
& \geq \sum_{j\in\sigma}\langle \Lambda_{j}Cf, \Lambda_{j}C'f\rangle
+
\sum_{j\in\sigma^{c}}  \sum_{k=1}^{\infty}|\alpha_{k}|^2|\beta_{kj}^{k}|^2 \left\langle  \Lambda_{j}(CC')^{\frac{1}{2}}(e_{k}), \Lambda_{j}(CC')^{\frac{1}{2}}(e_{k})\right\rangle\\
& \geq \sum_{j\in\sigma}\langle \Lambda_{j}Cf, \Lambda_{j}C'f\rangle
+
M\sum_{j\in\sigma^{c}}  \sum_{k=1}^{\infty}|\alpha_{k}|^2 \left\langle  \Lambda_{j}(CC')^{\frac{1}{2}}(e_{k}), \Lambda_{j}(CC')^{\frac{1}{2}}(e_{k})\right\rangle\\
& =\sum_{j\in\sigma}\langle \Lambda_{j}Cf, \Lambda_{j}C'f\rangle
+
M\sum_{j\in\sigma^{c}}  \left\langle  \sum_{k=1}^{\infty}\alpha_{k} \Lambda_{j}(CC')^{\frac{1}{2}}(e_{k}),  \sum_{k=1}^{\infty}\alpha_{k}\Lambda_{j}(CC')^{\frac{1}{2}}(e_{k})\right\rangle\\
& =\sum_{j\in\sigma}\langle \Lambda_{j}Cf, \Lambda_{j}C'f\rangle
+
M\sum_{j\in\sigma^{c}}  \left\langle   \Lambda_{j}(CC')^{\frac{1}{2}}\left(\sum_{k=1}^{\infty}\alpha_{k}(e_{k})\right),  \Lambda_{j}(CC')^{\frac{1}{2}}\left(\sum_{k=1}^{\infty}\alpha_{k}(e_{k})\right)\right\rangle\\
& =\sum_{j\in\sigma}\langle \Lambda_{j}Cf, \Lambda_{j}C'f\rangle
+
M\sum_{j\in\sigma^{c}}  \left\langle   \Lambda_{j}(CC')^{\frac{1}{2}}f,  \Lambda_{j}(CC')^{\frac{1}{2}}f\right\rangle\\
&\geq\min\{1,\ M\}\left( \sum_{j\in\sigma}\langle \Lambda_{j}Cf, \Lambda_{j}C'f\rangle
+
\sum_{j\in\sigma^{c}} \langle \Lambda_{j}Cf, \Lambda_{j}C'f\rangle  \right)\\
&= \min\{1,\ M\} \sum_{j=1}^\infty \langle \Lambda_{j}Cf, \Lambda_{j}C'f\rangle \\
&\geq A\min\{1,\ M\} \Vert K^{*}f\Vert^{2}.
\end{align*}
Hence $\Lambda$ and $\Omega$ are $K$-$g$-woven with universal $K$-$g$-frame bounds $(A\displaystyle \min\{1,\ M\})$ and $(B+\beta)$. $\square $

\begin{example}
Let $\{e_{k}\}_{k=1}^{\infty}$ be an orthonormal basis of a complex separable Hilbert space $\mathcal{H}$. Let $ K $ be the orthogonal projection from $\mathcal{H}$ onto $\overline{\mathrm{s}\mathrm{p}\mathrm{a}\mathrm{n}}\{e_{k}\}_{k=3}^{\infty}$ and $C, C'\in \mathcal{GL^{+}(H)}$ be such that
$$
C(e_{i}) =C'(e_{i}) = \left\{\begin{array}{ll}
e_{1}+ e_2, & i=1,\\
e_i, & \mathrm{otherwise}.
\end{array}\right.$$
 Define the sequences of bounded linear operators $\Lambda\equiv\{\Lambda_{j}\}_{j=1}^{\infty}$ and $\Omega\equiv\{\Omega_{j}\}_{j=1}^{\infty}$ on $\mathcal{H}$ as follows:
$$
\Lambda_{j}(e_{i})=\left\{\begin{array}{ll}
e_{i}, & i=j+2,\ j+3\\
0, & \mathrm{otherwise}.
\end{array}\right.
$$
and
$$
\Omega_{j}(e_{i})=\left\{\begin{array}{ll}
(1+\frac{1}{2^{j}})e_{j+2}, & i=j+2\\
e_{j+3}, & i=j+3\\
0, & \mathrm{otherwise}.
\end{array}\right.    (i,\ j\in \mathbb{N})
$$
Then, $\{\Lambda_{j}(CC')^{\frac{1}{2}}(e_{k})\}_{k=1}^{\infty}$ and $\{\Omega_{j}(CC')^{\frac{1}{2}}(e_{k})\}_{k=1}^{\infty}$ are orthogonal sets. Further, for all $f\in \mathcal{H}$ we have
$$
\Vert K ^{*}f\Vert^{2}\leq\sum_{j=1}^{\infty}\langle \Lambda_{j}Cf, \Lambda_{j}C'f\rangle \leq 2\Vert f\Vert^{2}.
$$
Thus, $\Lambda$ is a $(C, C')$-controlled $K$-$g$-frame for $\mathcal{H}$ with $(C, C')$-controlled $K$-$g$-frame bounds 1 and 2. Similarly, we can show that $\Omega$ is a  $(C, C')$-controlled $ K -$-frame with  $(C, C')$-controlled  $K$-$g$-frame bounds 1 and 5.

For $k\in \mathbb{N}$, let $\{\beta_{ij}^{k}\}_{i,j\in \mathbb{N}}$ be a sequence of scalars given by
$$
\beta_{ij}^{k}=\left\{\begin{array}{ll}
1, & i=k\neq j+2,\\
1+\frac{1}{2^{j}}, & i=k=j+2,\\
0, & \mathrm{o}\mathrm{t}\mathrm{h}\mathrm{e}\mathrm{r}\mathrm{w}\mathrm{i}\mathrm{s}\mathrm{e},
\end{array}\right.  (i,\ j\in \mathbb{N})
$$
Then, $\displaystyle \inf\{\beta_{kj}^{k}$ : $j\in \mathbb{N}\}=1>0$ and
$$
 \Omega_{j}(CC')^{\frac{1}{2}}(e_{k})=\sum_{i=1}^{\infty}\beta_{ij}^{k}\Lambda_{j}(CC')^{\frac{1}{2}}(e_i)
 $$
  for all $ j\in$ N. Hence by \ref{T3.6}, $\Lambda$ and $\Omega$ are  $(C, C')$-controlled $K$-$g$-woven.
\end{example}

An arbitrary $(C, C')$-controlled $ K $-$g$-Bessel sequence in $\mathcal{H}$ need not be a $(C, C')$-controlled $ K $-g- frame for $\mathcal{H}$. The following theorem gives sufficient conditions for $(C, C')$-controlled $ K $-$g$-Bessel sequences to constitute woven $(C, C')$-controlled $K$-$g$-frames for the underlying space.

\begin{theorem}
 Let $\Lambda\equiv\{\Lambda_{j}\}_{j=1}^{\infty}$  and $\Omega\equiv\{\Omega_{j}\}_{j=1}^{\infty}$  be  $(C, C')$-controlled  $ K $-$g$- Bessel sequences for $\mathcal{H}$  with respect to $\{\mathcal{H}_{j}\}_{j=1}^{\infty}$  such that for each $f\in \mathcal{H},$
$$  K (f)=\sum_{i\in \mathbb{N}}C\Lambda_{i}^{*}\Omega_{i}C'(f)\ \ \ \mbox{ and}\ \ \ C\Lambda_{i}^{*}\Omega_{i}C'=C\Omega_{i}^{*}\Lambda_{i}C'\ \  \ (i\in \mathbb{N}).$$
Then, $\Lambda$ and $\Omega$  are  $(C, C')$-controlled  $ K $-$g$-woven frames for $\mathcal{H}.$
\end{theorem}
\noindent\textbf{Proof.}
 Let $B_{1}$ and $B_{2}$ be $(C, C')$-controlled $ K $-$g$-Bessel bounds for $\Lambda$ and $\Omega$, respectively. Then, for any subset $\sigma$ of $\mathbb{N}$, we have
\begin{align*}
\Vert K ^{*}(f)\Vert^{4}
& =|\left\langle  K ^{*}(f),  K ^{*}(f)\right\rangle|^{2}\\
&=|\left\langle K  K ^{*}(f),\ f\right\rangle|^{2}\\
&=\left|\left\langle\sum_{i\in \mathbb{N}} C\Lambda_{i}^{*}\Omega_{i}C'( K ^{*}(f)),\ f\right\rangle\right|^{2}\\
&=\left|\left\langle\sum_{i\in \mathbb{N}} (CC')^{\frac{1}{2}}\Lambda_{i}^{*}\Omega_{i}(CC')^{\frac{1}{2}}( K ^{*}(f)),\ f\right\rangle\right|^{2}\\
&=\left|\left\langle\sum_{i\in\sigma} (CC')^{\frac{1}{2}}\Lambda_{i}^{*}\Omega_{i}(CC')^{\frac{1}{2}}( K ^{*}(f)) + \sum_{i\in\sigma^{c}} (CC')^{\frac{1}{2}}\Omega_{i}^{*}\Lambda_{i}(CC')^{\frac{1}{2}}( K ^{*}(f)) ,\ f\right\rangle\right|^{2}\\
&\leq 2\left|\left\langle\sum_{i\in\sigma} (CC')^{\frac{1}{2}}\Lambda_{i}^{*}\Omega_{i}(CC')^{\frac{1}{2}}( K ^{*}(f)) ,\ f\right\rangle\right|^{2}
+ 2\left|\left\langle\sum_{i\in\sigma^{c}} (CC')^{\frac{1}{2}}\Omega_{i}^{*}\Lambda_{i}(CC')^{\frac{1}{2}}( K ^{*}(f)) ,\ f\right\rangle\right|^{2}\\
&= 2\left| \sum_{i\in\sigma}\left\langle (CC')^{\frac{1}{2}}\Lambda_{i}^{*}\Omega_{i}(CC')^{\frac{1}{2}}( K ^{*}(f)) ,\ f\right\rangle\right|^{2}
+ 2\left|  \sum_{i\in\sigma^{c}} \left\langle (CC')^{\frac{1}{2}}\Omega_{i}^{*}\Lambda_{i}(CC')^{\frac{1}{2}}( K ^{*}(f)) ,\ f\right\rangle\right|^{2}\\
&= 2\left|\sum_{i\in\sigma} \left\langle \Omega_{i}(CC')^{\frac{1}{2}}( K ^{*}(f)) ,\ \Lambda_{i}(CC')^{\frac{1}{2}}f\right\rangle\right|^{2}
+ 2\left|\sum_{i\in\sigma^{c}}\left\langle  \Lambda_i(CC')^{\frac{1}{2}}( K ^{*}(f)) ,\ \Omega_{i}(CC')^{\frac{1}{2}}f\right\rangle\right|^{2}\\
&\leq  2\left|\sum_{i\in\sigma} \left\Vert \Omega_{i}(CC')^{\frac{1}{2}}( K ^{*}(f))\right\Vert  \left\Vert\Lambda_{i}(CC')^{\frac{1}{2}}f\right\Vert\right|^{2}
+ 2\left|\sum_{i\in\sigma^{c}}\left\Vert  \Lambda_{i}(CC')^{\frac{1}{2}}( K ^{*}(f))\right\Vert \left\Vert\Omega_i(CC')^{\frac{1}{2}}f\right\Vert\right|^{2}
\end{align*}
\begin{align*}
&\leq  2\sum_{i\in\sigma} \left\Vert \Omega_{i}(CC')^{\frac{1}{2}}( K ^{*}(f))\right\Vert^2  \sum_{i\in\sigma}\left\Vert\Lambda_{i}(CC')^{\frac{1}{2}}f\right\Vert^2\\
& \ \ + 2\sum_{i\in\sigma^{c}}\left\Vert \Lambda_{i}(CC')^{\frac{1}{2}}( K ^{*}(f))\right\Vert^2 \sum_{i\in\sigma^{c}} \left\Vert\Omega_{i}(CC')^{\frac{1}{2}}f\right\Vert^2\\
&= 2\left(\sum_{i\in\sigma} \left\langle \Omega_{i}(CC')^{\frac{1}{2}}( K ^{*}(f)) ,\ \Omega_{i}(CC')^{\frac{1}{2}}( K ^{*}(f))\right\rangle\right)
\left(\sum_{i\in\sigma} \left\langle \Lambda_{i}(CC')^{\frac{1}{2}}f ,\ \Lambda_{i}(CC')^{\frac{1}{2}}f\right\rangle\right)\\
&\ \ + 2\left(\sum_{i\in\sigma^{c}} \left\langle \Lambda_{i}(CC')^{\frac{1}{2}}( K ^{*}(f)) ,\ \Lambda_{i}(CC')^{\frac{1}{2}}( K ^{*}(f))\right\rangle\right)
\left(\sum_{i\in\sigma^{c}} \left\langle \Omega_{i}(CC')^{\frac{1}{2}}f ,\ \Omega_{i}(CC')^{\frac{1}{2}}f\right\rangle\right)\\
&= 2\left(\sum_{i\in\sigma} \left\langle  K ^{*}(f)  ,\ C\Omega_{i}^*\Omega_{i}C'( K ^{*}(f)) \right\rangle\right)
\left(\sum_{i\in\sigma} \left\langle f ,\ C\Lambda_{i}^*\Lambda_{i}C'f\right\rangle\right)\\
&\ \ + 2\left(\sum_{i\in\sigma^{c}} \left\langle  K ^{*}(f),\ C\Lambda_{i}^*\Lambda_{i}C'( K ^{*}(f))\right\rangle\right)
\left(\sum_{i\in\sigma^{c}} \left\langle f  ,\ C\Omega_{i}^*\Omega_{i}C'f\right\rangle\right)\\
&= 2\left(\sum_{i\in\sigma} \left\langle \Omega_{i}C( K ^{*}(f))  ,\ \Omega_{i}C'( K ^{*}(f)) \right\rangle\right)
\left(\sum_{i\in\sigma} \left\langle \Lambda_{i}Cf ,\ \Lambda_{i}C'f\right\rangle\right)\\
&\ \ + 2\left(\sum_{i\in\sigma^{c}} \left\langle \Lambda_{i}C( K ^{*}(f)),\ \Lambda_{i}C'( K ^{*}(f))\right\rangle\right)
\left(\sum_{i\in\sigma^{c}} \left\langle \Omega_{i}C f  ,\ \Omega_{i}C'f\right\rangle\right)\\
&\leq 2 B_2 \left\Vert K ^{*}(f)\right\Vert^2 \sum_{i\in\sigma} \left\langle \Lambda_{i}Cf ,\ \Lambda_{i}C'f\right\rangle + 2 B_1 \left\Vert K ^{*}(f)\right\Vert^2\sum_{i\in\sigma^{c}} \left\langle \Omega_{i}C f  ,\ \Omega_{i}C'f\right\rangle\\
&\leq 2 \max \{ B_1,B_2\} \left\Vert K ^{*}(f)\right\Vert^2 \left(\sum_{i\in\sigma} \left\langle \Lambda_{i}Cf ,\ \Lambda_{i}C'f\right\rangle + \sum_{i\in\sigma^{c}} \left\langle \Omega_{i}C f  ,\ \Omega_{i}C'f\right\rangle\right)\\
\end{align*}
Therefore, for all $f\in \mathcal{H}$, we have
$$ \frac{1}{2\max\{B_{1},B_{2}\}}\Vert K ^{*}(f)\Vert^{2}\leq\sum_{i\in\sigma} \left\langle \Lambda_{i}Cf ,\ \Lambda_{i}C'f\right\rangle + \sum_{i\in\sigma^{c}} \left\langle \Omega_{i}C f  ,\ \Omega_{i}C'f\right\rangle\leq(B_{1}+B_{2})\Vert f\Vert^{2}.$$
Hence $\Lambda$ and $\Omega$ are $(C, C')$-controlled $K$-$g$-woven. $\square $\\\\

The next theorem provides a necessary and sufficient condition for $(C, C')$-controlled $K$-g-woven frames which connects to ordinary weaving $K$-frames. More precisely, if frame bounds of frames associated with atomic spaces are positively confined, then $(C, C')$-controlled $K$-$g$-woven frames give ordinary weaving $K$-frames and vice-versa.
\begin{theorem}
Suppose that $\Lambda\equiv\{\Lambda_{j}\}_{j=1}^{\infty}$  and $\Omega\equiv\{\Omega_{j}\}_{j=1}^{\infty}$  are $(C, C')$-controlled  $K$-$g$-frames for $\mathcal{H}$  with respect to $\{\mathcal{H}_{j}\}_{j=1}^{\infty}$  and $\{\mathcal{W}_{j}\}_{j=1}^{\infty}$,  respectively. Let $\{f_{jk}\}_{k\in I_{j}\subset \mathbb{N}}$  and $\{g_{jk}\}_{k\in Q_{j}\subset \mathbb{N}}$ be frames for $\mathcal{H}_{j}$  and $\mathcal{W}_{j}$,  respectively $(j\in \mathbb{N})$  with frame bounds $\alpha, \beta$  and $\alpha', \beta'$,  respectively. Then the following conditions are equivalent.

(i) $\Lambda$  and $\Omega$  are $(C, C')$-controlled $K$-$g$-woven.

(ii) $\{(CC')^{\frac{1}{2}}\Lambda_{j}^{*}f_{jk}\}_{j\in \mathbb{N},k\in I_{j}}$ and $\{(CC')^{\frac{1}{2}}\Omega_{j}^{*}g_{jk}\}_{j\in \mathbb{N},k\in Q_{j}}$  are woven $K$-frames for $\mathcal{H}.$
\end{theorem}
\noindent\textbf{Proof.}  $(i)\Rightarrow(ii)$ : Let $A$ and $B$ be universal $(C, C')$-controlled $K$-$g$-frame bounds for $\Lambda$ and $\Omega.$ For any subset $\sigma$ of $\mathbb{N}$ and for any $f\in \mathcal{H}$, we compute
\begin{align*}
\sum_{j\in\sigma}\sum_{k\in I_{j}} & \left|\left\langle f, (CC')^{\frac{1}{2}}\Lambda_{j}^{*}f_{jk}\right\rangle\right|^{2}+\sum_{j\in\sigma^{c}}\sum_{k\in Q_{j}}\left|\left\langle f, (CC')^{\frac{1}{2}}\Omega_{j}^{*}g_{jk}\right\rangle\right|^{2}\\
&=\sum_{j\in\sigma}\sum_{k\in I_{j}}\left|\left\langle\Lambda_{j} (CC')^{\frac{1}{2}}(f), f_{jk}\right\rangle\right|^{2}+\sum_{j\in\sigma^{c}}\sum_{k\in Q_{j}}\left|\left\langle\Omega_{j} (CC')^{\frac{1}{2}}(f), g_{jk}\right\rangle\right|^{2}\\
&\leq\beta\sum_{j\in\sigma}\left\Vert\Lambda_{j} (CC')^{\frac{1}{2}}(f)\right\Vert^{2}+\beta'\sum_{j\in\sigma^{c}}\left\Vert\Omega_{j} (CC')^{\frac{1}{2}}(f)\right\Vert^{2}\\
&\leq\max\{\beta, \beta'\}\left(\sum_{j\in\sigma}\left\Vert\Lambda_{j} (CC')^{\frac{1}{2}}(f)\right\Vert^{2}+\sum_{j\in\sigma^{c}}\left\Vert\Omega_{j} (CC')^{\frac{1}{2}}(f)\right\Vert^{2}\right)\\
& = \max\{\beta, \beta'\}\left(\sum_{j\in\sigma}   \left\langle \Lambda_{j} (CC')^{\frac{1}{2}}(f), \Lambda_{j} (CC')^{\frac{1}{2}}(f)\right\rangle\right.\\
& \ \ \ \hspace{3cm}\left.
+\sum_{j\in\sigma^{c}}  \left\langle \Omega_{j} (CC')^{\frac{1}{2}}(f), \Omega_{j} (CC')^{\frac{1}{2}}(f)\right\rangle \right)\\
& = \max\{\beta, \beta'\}\left(\sum_{j\in\sigma}   \left\langle \Lambda_{j}C  f, \Lambda_{j}C'f\right\rangle
+\sum_{j\in\sigma^{c}}  \left\langle \Omega_{j}C  f, \Omega_{j}C' f\right\rangle \right)\\
& \leq\max\{\beta, \beta'\}B\Vert f\Vert^{2}.
\end{align*}
This gives a universal upper $K$-frame bound for the family $\{(CC')^{\frac{1}{2}}\Lambda_{j}^{*}f_{jk}\}_{j\in \mathbb{N},k\in I_{j}}$ and $\{(CC')^{\frac{1}{2}}\Omega_{j}^{*}g_{jk}\}_{j\in \mathbb{N},k\in Q_{j}}.$
For the universal lower $K$-frame bound, we compute
\begin{align*}
\sum_{j\in\sigma}\sum_{k\in I_{j}} & \left|\left\langle f, (CC')^{\frac{1}{2}}\Lambda_{j}^{*}f_{jk}\right\rangle\right|^{2}+\sum_{j\in\sigma^{c}}\sum_{k\in Q_{j}}\left|\left\langle f, (CC')^{\frac{1}{2}}\Omega_{j}^{*}g_{jk}\right\rangle\right|^{2}\\
&=\sum_{j\in\sigma}\sum_{k\in I_{j}}\left|\left\langle\Lambda_{j} (CC')^{\frac{1}{2}}(f), f_{jk}\right\rangle\right|^{2}+\sum_{j\in\sigma^{c}}\sum_{k\in Q_{j}}\left|\left\langle\Omega_{j} (CC')^{\frac{1}{2}}(f), g_{jk}\right\rangle\right|^{2}\\
&\geq\alpha\sum_{j\in\sigma}\left\Vert\Lambda_{j} (CC')^{\frac{1}{2}}(f)\right\Vert^{2}+\alpha'\sum_{j\in\sigma^{c}}\left\Vert\Omega_{j} (CC')^{\frac{1}{2}}(f)\right\Vert^{2}\\
&\geq\min\{\alpha, \alpha'\}\left(\sum_{j\in\sigma}\left\Vert\Lambda_{j} (CC')^{\frac{1}{2}}(f)\right\Vert^{2}+\sum_{j\in\sigma^{c}}\left\Vert\Omega_{j} (CC')^{\frac{1}{2}}(f)\right\Vert^{2}\right)\\
& = \min\{\alpha, \alpha'\}\left(\sum_{j\in\sigma}   \left\langle \Lambda_{j} (CC')^{\frac{1}{2}}(f), \Lambda_{j} (CC')^{\frac{1}{2}}(f)\right\rangle\right.\\
& \ \ \ \hspace{3cm}\left. +\sum_{j\in\sigma^{c}}  \left\langle \Omega_{j} (CC')^{\frac{1}{2}}(f), \Omega_{j} (CC')^{\frac{1}{2}}(f)\right\rangle \right)\\
& = \min\{\alpha, \alpha'\}\left(\sum_{j\in\sigma}   \left\langle \Lambda_{j}C  f, \Lambda_{j}C'f\right\rangle
+\sum_{j\in\sigma^{c}}  \left\langle \Omega_{j}C  f, \Omega_{j}C' f\right\rangle \right)\\
& \geq\min\{\beta, \beta'\}A\Vert K ^{*}f\Vert^{2}, f\in \mathcal{H}.
\end{align*}
Hence $\{ (CC')^{\frac{1}{2}}\Lambda_{j}^{*}f_{jk}\}_{j\in \mathbb{N},k\in I_{j}}$ and $\{ (CC')^{\frac{1}{2}}\Omega_{j}^{*}g_{jk}\}_{j\in \mathbb{N},k\in Q_{j}}$ are woven $K$-frames for $\mathcal{H}.$

$(ii)\Rightarrow(i)$ : Let $C$ and $D$ be universal $K$-frame bounds for $\{ (CC')^{\frac{1}{2}}\Lambda_{j}^{*}f_{jk}\}_{j\in \mathbb{N},k\in I_{j}}$ and $\{ (CC')^{\frac{1}{2}}\Omega_{j}^{*}g_{jk}\}_{j\in \mathbb{N},k\in Q_{j}}$. Then, for any subset $\sigma$ of $\mathbb{N}$, we compute
\begin{align*}
\sum_{j\in\sigma}  & \left\langle \Lambda_{j}C  f, \Lambda_{j}C'f\right\rangle
+\sum_{j\in\sigma^{c}}  \left\langle \Omega_{j}C  f, \Omega_{j}C' f\right\rangle\\
& = \sum_{j\in\sigma}   \left\langle \Lambda_{j} (CC')^{\frac{1}{2}}(f), \Lambda_{j} (CC')^{\frac{1}{2}}(f)\right\rangle +\sum_{j\in\sigma^{c}}  \left\langle \Omega_{j} (CC')^{\frac{1}{2}}(f), \Omega_{j} (CC')^{\frac{1}{2}}(f)\right\rangle \\
&= \sum_{j\in\sigma}\left\Vert\Lambda_{j} (CC')^{\frac{1}{2}}(f)\right\Vert^{2}+ \sum_{j\in\sigma^{c}}\left\Vert\Omega_{j} (CC')^{\frac{1}{2}}(f)\right\Vert^{2}\\
&\leq \frac{1}{\alpha}\sum_{j\in\sigma}\sum_{k\in I_{j}}\left|\left\langle\Lambda_{j} (CC')^{\frac{1}{2}}(f), f_{jk}\right\rangle\right|^{2}+\frac{1}{\alpha'}\sum_{j\in\sigma^{c}}\sum_{k\in Q_{j}}\left|\left\langle\Omega_{j} (CC')^{\frac{1}{2}}(f), g_{jk}\right\rangle\right|^{2}\\
& = \frac{1}{\alpha} \sum_{j\in\sigma}\sum_{k\in I_{j}}  \left|\left\langle f, (CC')^{\frac{1}{2}}\Lambda_{j}^{*}f_{jk}\right\rangle\right|^{2}+ \frac{1}{\alpha'} \sum_{j\in\sigma^{c}}\sum_{k\in Q_{j}}\left|\left\langle f, (CC')^{\frac{1}{2}}\Omega_{j}^{*}g_{jk}\right\rangle\right|^{2}\\
& \leq\max\left\{\frac{1}{\alpha}, \frac{1}{\alpha}\right\}D\Vert f\Vert^{2} \ \ \ \ \ \forall f\in \mathcal{H}.
\end{align*}
This gives a universal upper $(C, C')$-controlled  $K$-$g$-frame bound for the family $\Lambda$ and $\Omega$. Similarly, $C\displaystyle \min\{\frac{1}{\beta}, \displaystyle \frac{1}{\beta}\}$ will be a universal lower $K$-$g$-frame bound for $\Lambda$ and $\Omega$. The proof is complete. $\square $

The following theorem gives a sufficient condition for weaving $K$-$g$-frames in terms of positive operators associated with given $ K $-$g$ frames.

\begin{theorem}
Let $\Lambda\equiv\{\Lambda_{j}\}_{\}=1}^{\infty}$  and $\Omega\equiv\{\Omega_{j}\}_{j=1}^{\infty}$  be $K$-$g$-frames for $\mathcal{H}$ with respect to $\{\mathcal{H}_{j}\}_{j=1}^{\infty}$  and $\{\mathcal{W}_{j}\}_{j=1}^{\infty}$,  respectively. For any $\mathrm{J}\subseteq \mathbb{N}$,  suppose that the operator $U_{\mathrm{J}}:\mathcal{H}\rightarrow \mathcal{H}$  defined by
$$
U_{\mathrm{J}}(f)=\sum_{i\in \mathrm{J}}[\Omega_{i}^{*}\Omega_{i}(f)-\Lambda_{i}^{*}\Lambda_{i}(f)],\ f\in \mathcal{H},
$$
 is a positive linear operator. Then $\Lambda$  and $\Omega$  are $K$-$g$-woven.
\end{theorem}
\noindent\textbf{Proof.}
 Let $A, B$ and $\alpha, \beta$ be $K$-$g$-frame bounds for $\Lambda$ and $\Omega$, respectively. Then, for any subset $\sigma\subset \mathbb{N}$, we compute
\begin{align*}
A\Vert K ^{*}(f)\Vert^{2} & \leq \sum_{j\in \mathbb{N}} \left\langle \Lambda_{j}C  f, \Lambda_{j}C'f\right\rangle
\\
& =\sum_{j\in\sigma}  \left\langle \Lambda_{j}C  f, \Lambda_{j}C'f\right\rangle +\sum_{j\in\sigma^{c}}  \left\langle \Lambda_{j}C  f, \Lambda_{j}C'f\right\rangle\\
& =\sum_{j\in\sigma}  \left\langle \Lambda_{j}C  f, \Lambda_{j}C'f\right\rangle +  \left\langle   f, \sum_{j\in\sigma^{c}} C\Lambda_{j}^*\Lambda_{j}C'f\right\rangle\\
& =\sum_{j\in\sigma}  \left\langle \Lambda_{j}C  f, \Lambda_{j}C'f\right\rangle + \left\langle   f, \sum_{j\in\sigma^{c}}  C \Omega_{j}^*\Omega_{j}C'f-U_{\sigma^{c}}(f)\right\rangle
\\
& \leq\sum_{j\in\sigma}  \left\langle \Lambda_{j}C  f, \Lambda_{j}C'f\right\rangle + \left\langle   f, \sum_{j\in\sigma^{c}}  C \Omega_{j}^*\Omega_{j}C'f\right\rangle
\\
&= \sum_{j\in\sigma}\left\langle \Lambda_{j}C  f, \Lambda_{j}C'f\right\rangle  +\sum_{j\in\sigma^{c}}  \left\langle \Omega_{j}C  f, \Omega_{j}C'f\right\rangle  \\
& \leq(B+\beta)\Vert f\Vert^{2} \ \ \ \ \ \forall f\in \mathcal{H}.
\end{align*}
Hence $\Lambda$ and $\Omega$ are $K$-$g$-woven with universal $K$-$g$-frame bounds $A$ and $(B+ \beta)$ . $\square $

\par
\textbf{Acknowledgments:} The authors would like to thank referee(s) for valuable comments and suggestions.

\textbf{References}









\end{document}